\documentclass[reqno]{amsart}
\usepackage{latexsym}
\usepackage{amssymb,amsthm,amsmath}
\usepackage[notref,notcite]{showkeys}

\theoremstyle{plain}
\newtheorem{theorem}{Theorem}
\newtheorem{lemma}[theorem]{Lemma}

\theoremstyle{definition}

\theoremstyle{remark}
\newtheorem{remark}[theorem]{Remark}

\def\d#1{{#1\kern-0.4em\char"16\kern-0.1em}}
\def\D#1{{\raise0.2ex\hbox{-}\kern-0.4em #1}}
\newcounter{zd}

\newcounter{zdr}[subsection]

\newcommand{\eps}{\varepsilon}
\DeclareMathOperator{\Pim }{\Pi}

\def\R{I\!\!R}

\def\N{I\!\!N}

\def\F{{\cal F}}

\def\cal{\mathcal}

\begin{document}

\title[On the first commutation lemma]
{On a variant of Tartar's first commutation lemma}
\author{ D.~ Mitrovi\'c}
\address{ Darko Mitrovic, University of Montenegro, Faculty of Mathematics, Cetinjski put bb, 81000 Podgorica, Montenegro}
 \email{  matematika@t-com.me}
 \date{}

\begin{abstract}
We prove a variant of Tartar's first commutation lemma involving
multiplier operators with symbols not necessarily defined on a
manifold of codimension one.
\end{abstract}


\maketitle
\section{Introduction}

At the beginning of 90's, L.Tartar \cite{Tar} and P.Gerard
\cite{Ger} independently introduced the $H$-measures (microlocal
defect measures). The $H$-measures appeared to be very powerful tool
in many fields of mathematics and physics (see randomly chosen
\cite{?7, ?5, ?2, ?9, ?3, JMSpa, ?4, ?1}). They are given by the
following theorem:

\begin{theorem}\cite{Tar}
\label{tbasic1} If $(u_n)=((u_n^1,\dots, u_n^r))$, is a sequence in
$L^2(\R^d;\R^r)$ such that $u_n\rightharpoonup 0$ in
$L^2(\R^d;\R^r)$, then there exists its subsequence $(u_{n'})$ and a
positive definite matrix of complex Radon measures
$\mu=\{\mu^{ij}\}_{i,j=1,\dots,d}$ on $\R^d\times S^{d-1}$ such that
for all $\varphi_1,\varphi_2\in C_0(\R^d)$ and $\psi\in C(S^{d-1})$:
\begin{equation}
\label{basic1}
\begin{split}
\lim\limits_{n'\to \infty}\int_{\R^d}&(\varphi_1
u^i_{n'})(x)\overline{{\cal A}_\psi(\varphi_2
u^j_{n'})(x)}dx=\langle\mu^{ij},\varphi_1\overline{\varphi_2}\psi
\rangle\\&= \int_{\R^d\times
S^{d-1}}\varphi_1(x)\overline{\varphi_2(x)}\psi(\xi)d\mu^{ij}(x,\xi),
\ \ i,j=1,\dots,d,
\end{split}
\end{equation}where ${\cal A}_\psi$ is a multiplier operator with the symbol $\psi\in C^\kappa(S^{d-1})$.
\end{theorem}

The complex matrix Radon measure $\{\mu^{ij}\}_{i,j=1,\dots r}$
defined in the previous theorem we call the {\em $H$-measure}
corresponding to the subsequence $(u_{n'})\in L^2(\R^d;\R^r)$.

The $H$-measures describe a loss of strong $L^2$ precompactness for
the corresponding sequence $(u_n)\in L^2(\R^d;\R^r)$. In order to
describe loss of $L^1_{loc}$ precompactness for a sequence $(u_n)\in
L_{loc}^p(\R^d)$, $p>1$, we have extended the notion of the
$H$-measures in \cite{MPB} as follows.

\begin{theorem} \cite{MPB}
\label{tbasic2} Assume that $(u_n)$, is a sequence in
$L_{loc}^{p}({\R}^{d})$, $p>1$, such that $u_n\rightharpoonup 0$,
$n\to \infty$, in $L_{loc}^{p}(\R^{d})$, $\beta>0$. Assume that
$(v_n)$ is a bounded sequence in $L^\infty(\R^{d})$.

Then, there exist subsequences $(u_{n'})$ and $(v_{n'})$ of the
sequences $(u_{n})$ and $(v_n)$, respectively, such that there
exists a complex valued distribution $\mu\in {\cal D}'(\R^{d'}
\times S^{d-1})$, such that for every $\varphi_1, \varphi_2\in
C_0(\R^{d})$ and $\psi\in C^{\kappa}(S^{d-1})$, $\kappa>d/2$,
$\kappa\in \N$:
\begin{equation}
\label{basic2}
\begin{split}
\lim\limits_{n'\to \infty}\int_{\R^{d}}(\varphi_1 u_{n'})(x)&{{\cal
A}_{\psi}(\varphi_2 v_{n'})(x)}dx=\langle\mu ,\varphi_1\varphi_2\psi
\rangle,
\end{split}
\end{equation}where ${\cal A}_{\psi}:L^p(\R^d)\to L^p(\R^d)$ is a multiplier operator with the symbol $\psi\in C^\kappa(S^{d-1})$.
\end{theorem}

The first commutation lemma  was one of the key point in the proof
of Theorem \ref{tbasic1} as well as Theorem \ref{tbasic2} (with a
simple modification; see \cite[Lemma 14]{MPB}). It is stated as
follows:

\begin{lemma}\cite[Lemma 14]{MPB}
\label{scl} (First commutation lemma)  Let $a\in C(S^{d-1})$ and
$b\in C_0(\R^d)$. Let ${\cal A}$ be a multiplier operator with the
symbol $a$, and $B$ be an operator of multiplication given by the
formulae:
\begin{align*}
&\F({\cal A}u)(\xi)=a(\frac{\xi}{|\xi|})\F(u)(\xi) \ \ a.e. \ \ \xi\in \R^d,\\
 &Bu(x)=b(x)u(x) \ \ a.e. \ \ x\in \R^d,
\end{align*} where ${cal F}$ is the Fourier transform.
 Then $C={\cal A}B-B{\cal A}$ is a
compact operator from $L^p(\R^d)$ into $L^p(\R^d)$, $p>1$.
\end{lemma}

The proof of the lemma heavily relies on the fact that the function
$a$ is actually defined on the unit sphere. Recently, two new
variants of the $H$-measures were introduced -- the parabolic
$H$-measures \cite{Ant2} and ultra-parabolic $H$-measures
\cite{JMSpa}. In both cases a variant of the first commutation lemma
is needed, and in both cases its proof is based on the fact that a
symbol $a$ of appropriate multiplier ${\cal A}$ is defined on a
smooth, bounded, simply connected manifold of codimension one.

In order to motivate our variant of the first commutation lemma,
notice that from the proof of Theorem \ref{tbasic1} (see
\cite{Tar}), it follows that we need to "commute" ${\cal
A}(\varphi_2 u_n)$ by $\varphi_2 {\cal A}(u_n)$, where $(u_n)$ is
the sequence bounded in $L^2(\Omega)$, which was exactly done in
Lemma \ref{scl}. Similarly, in Theorem \ref{tbasic2} (see
\cite{MPB}), we need to "commute" ${\cal A}(\varphi_2 v_n)$ by
$\varphi_2 {\cal A}(v_n \chi_{{\rm supp \varphi_2}})$, where
$\chi_{V}$ is the characteristic function of the set $V$, and
$(v_n)$ is the sequence bounded in $L^\infty(\Omega)$. Therefore, it
is enough to prove that the commutator $C$ is compact operator from
$L^\infty_0(\Omega)$ into $L^p_{loc}(\Omega)$, $p>1$. We state:

\begin{lemma}
\label{scl'} Let $a\in C^{\kappa}(\R^d)$, $\kappa=\lfloor
d/2\rfloor+1$, and $b\in C_0(\R^d)$. Suppose that for some constant
$k>0$ and for any real number $r>0$
\begin{equation}
\label{cond_1'} \int_{\frac{r}{2}\leq \|\xi\|\leq
r}|D_{\xi}^{\alpha}a(\xi)|^2d\xi\leq k^2r^{d-2n(\alpha)}
\end{equation} holds for every $\alpha=(\alpha_1,\dots, \alpha_d)\in \N_0^d$ satisfying $n(\alpha)=\sum\limits_{i=1}^d\alpha_i\leq \kappa$.

Let ${\cal A}$ be a multiplier operator with the symbol $a$, and $B$
be an operator of multiplication given by the formulae:
\begin{align}
\label{oper_1}
&\F({\cal A}u)(\xi)=a(\xi)\F(u)(\xi) \ \ a.e. \ \ \xi\in \R^d,\\
\label{oper_2} &Bu(x)=b(x)u(x) \ \ a.e. \ \ x\in \R^d,
\end{align} where ${\cal F}$ is the Fourier transform.
 Then $C={\cal A}B-B{\cal A}$ is a compact operator from
$L_0^\infty(\R^d)$ into $L_{loc}^{p_0}(\R^d)$ for every
$1<p_0<\infty$.
\end{lemma}

\begin{remark}
We hope that the lemma could serve for defining variants of the
$H$-measures adapted to equations which change type (such as
non-strictly parabolic equations).
\end{remark}

\begin{proof}
On the first step notice that $a$ satisfies conditions of the
H\"{o}rmander-Mikhlin theorem (see \cite{Oki, Gra}). Therefore, for
every $p>1$ there exists a constant $k_p$ such that $\|{\cal
A}\|_{L^p\to L^p}\leq k_p$. Thus,
$$
\|C\|\leq 2 k_p \|b\|_{L^\infty(\R^d)}.
$$

Then, notice that we can assume $b\in C^1_0(\R^d)$. Indeed, if we
assume merely $b\in C_0(\R^d)$ then we can uniformly approach the
function $b$ by a sequence $(b_n)\in C^1_0(\R^d)$. The corresponding
sequence of commutators $C_n={\cal A} B_n-B_n {\cal A}$, where
$B_n(u)=b_n u$, converges in norm toward $C$. So, if we prove that
$C_n$ are compact for each $n$, the same will hold for $C$ as well.

Then, fix a real non-negative function $\omega$ with a compact
support and total mass one. Take the characteristic function
$\chi_{B(0,2)}$ of the ball $B(0,2)\subset \R^d$ and denote:
$$
\chi(x)=\chi_{B(0,2)}\star
\frac{1}{\eps^d}\Pim\limits_{i=1}^d\omega(\frac{x_i}{\eps})
$$for an $\eps>0$ small enough so that we have
$\chi(x)=1$ for $x\in B(0,1)$, and $(1-\chi)\equiv 1$ out of the
ball $B(0,3)$.

Next, notice that $ {\cal A}={\cal A}_{a \chi}+{\cal A}_{a
(1-\chi)}$, where ${\cal A}_{a \chi}$ is a multiplier operator with
the symbol $a \chi$, and ${\cal A}_{a (1-\chi)}$ is a multiplier
operator with the symbol $a (1-\chi)$. Accordingly,
\begin{align*}
C={\cal A}B-B{\cal A}&={\cal A}_{a \chi}B-B{\cal A}_{a \chi}
\\&+{\cal A}_{a (1-\chi)} B-B {\cal A}_{a (1-\chi)}=
C_{a \chi}+C_{a (1-\chi)},
\end{align*}where
$C_{a \chi}= {\cal A}_{a \chi}B-B{\cal A}_{a \chi}$ and $C_{a
(1-\chi)}={\cal A}_{a (1-\chi)} B-B {\cal A}_{a (1-\chi)}$.

First, consider the commutator $C_{a \chi}$.  Notice that since $a
\chi_{B(0,2)}$ has a compact support, the multiplier ${\cal A}_{a
\chi}$ is actually the convolution operator with the kernel
$\psi=\bar{\cal F}(a\chi)\in L^2(\R^d)$, where $\bar{\cal F}$ is the
inverse of the Fourier transform ${\cal F}$:
\begin{equation}
\label{oper_1'} {\cal A}_{a \chi}(u)=\psi\star u, \ \ u\in
L^\infty_0(\R^d).
\end{equation} Therefore, we can state that
\begin{equation}
\label{comm_1}
\begin{split}
&\qquad\qquad C_{a \chi}u(x)=\int_{\R^d}\left(b(x)-b(y)\right)\psi(x-y)u(y)dy,\\
&\text{is a compact operator from $L_0^\infty(\R^d)$ into
$L_{loc}^{p_0}(\R^d)$, $p_0\geq 2$.}
\end{split}
\end{equation}

Indeed, take an arbitrary bounded sequence $(u_n)\in
L_0^\infty(\R^d)$ such that $u_n\rightharpoonup 0$ weak-$\star$ in
$L^\infty(\R^d)$ and ${\rm supp}u_n\subset \hat{V}\subset\subset
\R^d$, for a relatively compact set $\hat{V}$. In order to prove
that $C_{a \chi}$ is compact, it is enough to prove that $C_{a
\chi}u_n$ strongly converges to zero in $L_{loc}^{p_0}(\R^d)$,
$p_0>1$.

Since $\psi\in L^2(\R^d)$ we also have $\psi\in L^1_{loc}(\R^d)$.
Thus, it holds for every fixed $x\in \R^d$
\begin{equation}
\label{comm_1-c1} C_{a
\chi}u_n(x)=\int_{\R^d}\left(b(x)-b(y)\right)\psi(x-y)u_n(y)dy\to 0,
\ \ n\to 0.
\end{equation} Next, since the sequence $(u_n)$ has compact support, we also have:
\begin{equation}
\label{comm_1-c2} |C_{a \chi}u_n(x)|\leq \hat{C},
\end{equation}for a constant $\hat{C}$ depending on the support of the sequence $(u_n)$ as well as $L^2$ norm of the kernel $\psi$.

Combining \eqref{comm_1-c1} and \eqref{comm_1-c2} with the Lebesgue
dominated convergence theorem, we see that for an arbitrary
relatively compact $V\subset\subset \R^d$ and every $p_0>0$, it
holds:
\begin{equation}
\label{comm_1-c3} \int_V|C_{a \chi}u_n(x)|^{p_0}dx\to 0, \ \ n\to
\infty,
\end{equation}proving \eqref{comm_1}.

In order to prove that $C_{a (1-\chi)}$ is compact, we need more
subtle arguments basically involving techniques from the proof of
the Hormander-Mikhlin theorem from e.g. \cite{Oki}.

So, let $\Theta$ be a non-negative infinitely differentiable
function supported by $\{\xi\in R^n:\frac{1}{2}\leq\|\xi\|\leq2\}$.
Also let $ \theta(\xi)=\Theta(\xi)\bigg/
\sum_{j=-\infty}^{\infty}\Theta(2^{-j}\xi)$, $\theta(0):=0$. Then,
$\theta$ is non-negative, it is supported by $\{\xi\in
R^n:\frac{1}{2}\leq \|\xi\|\leq2\}$, it is infinitely
differentiable, and is such that if $\xi\neq 0$, then $
\sum_{j=-\infty}^{\infty}\theta(2^{-j}\xi)=1. $

Now, let $a_j(\xi)=a(\xi)(1-\chi(\xi))\theta(2^{-j}\xi)$, $j>0$.
Then, ${a}_j$ has support in the set $ \{\xi \in R^n:
2^{j-1}\leq\|\xi\|\leq 2^{j+1}\}$, $j>0,$ and $
a(\xi)(1-\chi(\xi))=\sum_{j=0}^\infty {a}_j(\xi)$. Furthermore, it
holds for any $p\in \{1,\dots,d\}$, $\alpha_p\in \N_0$ that
$$
\frac{\partial^{\alpha_p}}{\partial^{\alpha_p}\xi_p}(a(\xi)(1-\chi(\xi))\theta(2^{-j}\xi))=\sum_{l=0}^{\alpha_p}{\alpha_p\choose
l}\frac{\partial^l a(\xi)}{\partial^l \xi_p}
\frac{\partial^{\alpha_p-l}(1-\chi(\xi))\theta(2^{-j}\xi)}{\partial^{\alpha_p-l}\xi},
$$ so that, with suitable bounded functions $a_{\beta\gamma}$,
$\beta+\gamma=\alpha$, $\alpha\in \N_0^d$, we have $
D_\xi^{\alpha}{\phi}_j(\xi)=\sum_{\beta+\gamma=
\alpha}a_{\beta\gamma}2^{-jn(\beta)}D_\xi^{\gamma}{a}(\xi)$.

From here, on applying hypothesis \eqref{cond_1'} with $r=2^j$, it
follows from the Minkowski inequality that
\begin{align}
\label{corr_ns1} \int_{R^d}|D^{\alpha}_\xi{a}_j(\xi)|^2d\xi  &\leq
p^2_0\sum_{\beta+\gamma=\alpha}a_{\beta\gamma}2^{-jn(\beta)}\int_{2^{j-1}
\leq\|\xi\|\leq 2^{j+1}}|D^{\gamma}_\xi{a}(\xi)|^2d\xi
\\&\leq C k^2
2^{j(d-2n(\alpha))}, \nonumber
\end{align}
 where $C$ is a constant independent on $k$.

Denote by $\bar{a}_j=\bar{\F}(a_j)(x)$, $x\in \R^d$, the inverse
Fourier transform of the function $a_j$. From \eqref{corr_ns1}, the
Cauchy-Schwartz inequality, Plancharel's theorem and the well known
properties of the Fourier transform, for every $s>0$ it holds (see
also the proof of \cite[Theorem 7.5.13]{Oki}):
\begin{align}
\label{aux_5}
\int_{\|x\|>s}|\bar{a}_j(x)|dx&\leq\bigg(\int_{\|x\|\geq s}
\|x\|^{-2\kappa}dx\bigg)^{{1}/{2}} \bigg(\int_{\|x\|\geq s}
\|x\|^{2\kappa}|\bar{a}_j(x)|^2d\xi\bigg)^{{1}/{2}}
\\&\leq \left(\frac{2\pi^{d-1}s^{d-2\kappa}}{2\kappa-d} \right)^{1/2} \left(d^\kappa \sum\limits_{i=1}^d
\int_{\R^d}|x_i|^{2\kappa}|\bar{a}_j(x)|^2dx\right)^{1/2}
\nonumber\\
&= \left(\frac{2\pi^{d-1}s^{d-2\kappa}}{2\kappa-d} \right)^{1/2}
\left(d^\kappa \sum\limits_{i=1}^d
\int_{\R^d}|D^{\kappa}_{\xi_i}{a}_j(\xi)|^2d\xi\right)^{1/2}
\nonumber\\
&\leq C_1 k (2^j s)^{(\frac{1}{2}d-\kappa)}, \nonumber
\end{align}where $C_1$ depends only on the functions $\theta$ and $\chi$.

%From here it follows
%\begin{equation}
%\label{sep128}
%\int_{\|x\|<s}\|x\||\bar{a}_j(x)|dx \to 0, \ \ s\to 0,
%\end{equation}since $\frac{1}{2}d-\kappa\geq -1$.

Next, consider $ \bar{A}_n(x)=\sum_{j=0}^n\bar{a}_j(x)$, $x\in
\R^d$. For an arbitrary fixed $s>0$, the series
$\sum_{j=0}^n\bar{a}_j(x)$ is absolutely convergent in
$L^1(\R^d\backslash B(0,s))$. Indeed,
\begin{align}
\label{sep138} \|\bar{A}_n(x)\|_{L^1(\R^d\backslash B(0,s))}&\leq
\sum_{j=0}^n\|\bar{a}_j\|_{L^1(\R^d\backslash B(0,s))}\\&\leq C_1 k
s^{(\frac{1}{2}d-\kappa)}\sum\limits_{j=0}^n
2^{j(\frac{1}{2}d-\kappa)}\leq C_3 <\infty, \nonumber
\end{align}for a constant $C_3>0$, since $\frac{1}{2}d-\kappa < 0$.

Thus, for every $s>0$ there exists $\bar{A}_s\in L^1(\R^d\backslash
B(0,s))$ such that
\begin{equation}
\label{sep148} \sum_{j=0}^\infty\bar{a}_j(x)=\bar{A}_s(x), \ \ x\in
\R^d\backslash B(0,s).
\end{equation}

Furthermore, for an odd $d$, from \eqref{aux_5} we have
\begin{equation}
\label{sep128} \int_{\|x\|<s}\|x\|\cdot|\bar{A}_n(x)|dx \to 0, \ \
s\to 0,
\end{equation}while for an even $d$ conclusion \eqref{sep128}
follows from:
\begin{align*}
\int_{\|x\|<s}\|x\|\cdot|\bar{A}_n(x)|dx&\leq
\left(\int_{\|x\|<s}\|x\|^{1-2\kappa}
\right)^{1/2}\left(\int_{\|x\|<s}\|x\|^{2\kappa-2}|\bar{a}_j(x)|^2dx\right)^{1/4}\times\\&\times
\left(\int_{\|x\|<s}\|x\|^{2\kappa}|\bar{a}_j(x)|^2dx\right)^{1/4}
\leq C_4 s^{1/2}.
\end{align*}

Now, take the convolution operator:
$$
{\cal A}_n(u)=\bar{A}_n\star u, \ \ u\in L^\infty_0(\R^d).
$$ and consider the commutator
$ C_n={\cal A}_n B-B {\cal A}_n$. It holds:
\begin{equation*}
C_n(u)(x)=-\int_{\R^d}\bar{A}_n(x-y)(b(x)-b(y))u(y)dy.
\end{equation*}Given a fixed $s>0$, rewrite $C_n(u)$ in the following way:
\begin{align*}
&C_n(u)(x)=\int_{\R^d}\bar{A}_n(x-y)(b(x)-b(y))u(y)dy\\
&=\int_{\|x-y\|>s}\bar{A}_n(x-y)(b(x)-b(y))u(y)dy+\int_{\|x-y\|\leq
s}\bar{A}_n(x-y)(b(x)-b(y))u(y)dy.
\end{align*} From here, combining $b\in C^1_0(\R^d)$ with \eqref{sep148} and \eqref{sep128}, we conclude for an arbitrary relatively compact
$V\subset\subset \R^d$:
\begin{equation}
\limsup\limits_{n\to \infty} \| C_n(u)(x)\|_{L^{p_0}(V)}\leq
\|\int_{\|x-y\|>s}
\bar{A}_s(x-y)(b(x)-b(y))u(y)dy\|_{L^{p_0}(V)}+o_s(1),
\end{equation}where
\begin{align*}
o_s(1)&=\int_{\|x-y\|\leq s}\bar{A}_n(x-y)(b(x)-b(y))u(y)dy\\&
=\int_{\|x-y\|\leq
s}\|x-y\|\bar{A}_n(x-y)\frac{(b(x)-b(y))}{\|x-y\|}u(y)dy
\to^{\eqref{sep128}} 0, \ \ \eps\to 0.
\end{align*} Furthermore, arguing as for \eqref{comm_1}, we infer that the operator
\begin{equation}
\label{nk1} u\mapsto \int_{\|x-y\|>s}\bar{A}_s(x-y)(b(x)-b(y))u(y)dy
\end{equation}is a compact operator from $L^\infty_0(\R^d)$ to
$L^{p_0}_{loc}(\R^d)$ for an arbitrary $p_0>1$.

Next, notice that we have for any $u\in L^2$:
\begin{align*}
&\|(C_n-C_{\phi (1-\chi)})(u)\|_{L^2}\\&\leq \|({\cal A}_n-{\cal
A}_{a (1-\chi)}(bu)\|_{L^2}+\|b\|_{\infty} \|({\cal A}_n-{\cal
A}_{\phi
(1-\chi)}(u)\|_{L^2}\\
&=\|(A_n-a(1-\chi))\F(bu)\|_{L^2}+\|b\|_{\infty}\|(A_n-a(1-\chi))\F(u)\|_{L^2}
\to 0, \ \ n\to\infty,
\end{align*}
and from here and the Fatou lemma:
\begin{align*}
&\|C_{a (1-\chi)}(u)\|_{L^{p_0}(V)}\leq
\limsup\limits_{n\to\infty}\bigg(\|\int_{\|x-y\|\geq
s}\bar{A}_n(x-y)(b(x)-b(y))u(y)dy\\&+\int_{\|x-y\|<
s}|x-y|\bar{A}_n(x-y)\frac{(b(x)-b(y))}{|x-y|}u(y)dy\|_{L^{p_0}(V)}
\bigg)\\&\leq \|\int_{\|x-y\|\geq
s}\bar{A}_s(x-y)(b(x)-b(y))u(y)dy\|_{L^{p_0}(V)}+ o_s(1),
\end{align*} where $o_s(1)$ denotes a quantity tending to zero as $s\to 0$, and appears here due to
\eqref{sep128}. From here and \eqref{nk1}, it follows that
$C_{a(1-\chi)}$ is a compact operator since it can be estimated by a
sum of a compact operator, and an operator bounded by an arbitrary
small constant.

Thus, we see that $C$ can be represented as the sum of two compact
operators $C_{a \chi}$ and $C_{a(1-\chi)}$, which means that $C$ is
a compact operator itself.

This concludes the proof.
\end{proof}


\begin{thebibliography}{99}
\bibitem{?7} G.~Alberti, S.~M\"{u}ller, {A new approach to
variational problems with multiple scales}, Comm. Pure. Appl. Math.
54(2001), 761-825

\bibitem{Ant2} N.~Antonic, M.~Lazar {\em $H$-measures and variants applied to parbolic equations},
J. Math. Anal. Appl., 343 (2008), 207-225

%\bibitem{Duoa} J.~Duoandikoetxea, {\em Fourier Analysis}, Graduate
%Studies in Mathematics, Volume 29, American Mathematical Society
%(2000)

%\bibitem{?6} W.~Bao, S.~Jin, P.~Markowich, {\em Numerical study of
%time-splitting spectral discretizations of nonlinear Schr\"{o}dinger
%equations in the semiclassical regimes}, SIAM J. Sci. Comput.
%25(2003), 27-64



\bibitem{Ger} Gerard, P.,
{\em Microlocal Defect Measures}, Comm. Partial Differential
Equations 16(1991), 1761--1794.

%\bibitem{HKM} H.~Holden, K.~Karlsen, D.~Mitrovic, {\em Zero
%diffusion dispersion limits for a
% scalar conservation law with
%discontinuous flux function}, preprint

%\bibitem{D}
%D.~Mitrovic, {\em On a strong precompactness of velocity averages for a heterogenous transport equation with rough coefficients}, preprint

%\bibitem{MP} D.~Mitrovic, S.~Pilipovic, {\em Approximations of linear Dirichlet problems
%with singularities}, J. Math. Anal. Appl. 313 (2006) 98-119.

\bibitem{Gra}L.~Grafakos, {\em Classical and Modern Fourier
Analysis}, Paerson Education, Inc., Upper Saddle River, New Jersey
07458 (2004)

\bibitem{?5} Y.~Brenier, {\em Hydrodynamic structure of the
augmented Born-Infeld equations}, Arch.Ration.Mech.Anal. 172(2004),
65-91

\bibitem{?2} R.~Lewandowski, {\em Vorticities in a LES model for 3D
periodic turbulent flows}, J.Math.Fluid Mech 8(2006), 398-442

\bibitem{?9} G.~Metivier, S.~Schochet, {\em Trilinear resonant
interactions of semilinear hyperbolic waves} Duke Math. J. 95
(1998), 241-304

\bibitem{MPB} D.~Mitrovic, N.~Antonic, S.~Pilipovic, V.~Bojkovic, {\em $H$-distributions -- an extension of the
$H$-measures}, arXiv:0907.1373v1

\bibitem{?3} A.~Mielke, {\em Macroscopic behaviour of microscopic
oscillations in harmonic lattices via Wigner-Husimi transforms},
Arch.Ration.Mech.Anal 181(2006), 401-448


\bibitem{MA} A.~Moulahi,
{\em Stabilisation interne d'ondes electromagnetiques dans un
domaine exterieur}, J. Math. Pures Appl. (9) 88 (2007), 431--453

\bibitem{Oki} G.~O.~Okikiolu,
{\em Aspects of the Theory of Bounded Integral operators in
$L^p$-Spaces}, Academic Press, London and new York, 1971

\bibitem{JMSpa}
E.Yu.~Panov, Ultra-parabolic equations with rough coefficients.
Entropy solutions and strong pre-compactness property, Journal of
Mathematical Sciences, 159:2(2009) 180--228.

\bibitem{?4} P.~Pedregal, {\em $\Gamma$ convergence through Young
measures}, SIAM J. Math. ANal. 36 (2004), 423-440

\bibitem{Tar} L.~Tartar, {\em H-measures, a new approach for studying homogenisation,
oscillation and concentration effects in PDEs}, Proc. Roy. Soc.
Edinburgh. Sect. A 115:3--4 (1990) 193-230

\bibitem{TL} L.~Thevenot, {\em An optimality condition for the assembly distribution in a
nuclear reactor}, Math. Models Methods Appl. Sci. 15(2005),
407--435.

\bibitem{?1}L.~Ying, E.~Candes, {\em The phase flow method},
J.Comput. Phys., 220(2006), 184-215

\end{thebibliography}
\end{document}